\documentclass[10pt,a4paper]{article}
\usepackage{graphics}
\usepackage{amsmath,amssymb,amsthm,graphics,graphicx,epic,eepic,multicol,ascmac}
\usepackage[mathscr]{euscript}
\usepackage[all]{xy}
\usepackage{fancybox}
\usepackage{url}
\usepackage{ytableau}
\usepackage{here}
\usepackage[toc,page]{appendix}

\ytableausetup{centertableaux,mathmode,boxsize=1em}
\def\young(#1){\ytableaushort{#1}}
\def\yng(#1){\ydiagram{#1}}

\numberwithin{equation}{section}
\theoremstyle{theorem}
\newtheorem{thm}{Theorem}[section]

\newtheorem{rem}[thm]{Remark}
\newtheorem{cor}[thm]{Corollary}

\theoremstyle{definition}
\newtheorem{defn}[thm]{Definition}
\newtheorem{ex}[thm]{Example}

\def\wht(#1){\widehat{\ #1\ }}

\newcommand{\ch}{\mathrm{ch}}

\newcommand{\lbr}{\begin{bmatrix}}
\newcommand{\rbr}{\end{bmatrix}}

\def\beneme{\begin{enumerate}}
\def\beq{\begin{equation}}
\def\beqn{\begin{eqnarray}}
\def\beqnn{\begin{eqnarray*}}

\def\bfii0{{\bf i_0}}

\def\bbra#1,#2,#3{\left\{\begin{array}{c}\hspace{-5pt}
#1;#2\\ \hspace{-5pt}#3\end{array}\hspace{-5pt}\right\}}

\def\ci(#1,#2){c_{#1}^{(#2)}}
\def\Ci(#1,#2){C_{#1}^{(#2)}}
\def\mpp(#1,#2,#3){#1^{(#2)}_{#3}}
\def\bCi(#1,#2){\ovl C_{#1}^{(#2)}}
\def\ch(#1,#2){c_{#2,#1}^{-h_{#1}}}
\def\cc(#1,#2){c_{#2,#1}}

\def\di(#1,#2){D_{#1}^{(#2)}}
\def\dbi(#1,#2){\ovl D_{#1}^{(#2)}}

\def\eneme{\end{enumerate}}

\def\eeq{\end{equation}}
\def\eeqn{\end{eqnarray}}
\def\eeqnn{\end{eqnarray*}}

\def\gau#1,#2{\left[\begin{array}{c}\hspace{-5pt}#1\\
\hspace{-5pt}#2\end{array}\hspace{-5pt}\right]}

\def\ji(#1,#2){j_{#1}^{(#2)}}

\def\lan{\langle}

\def\nd{\noindent}

\def\ovl{\overline}

\def\qed{\hfill\framebox[2mm]{}}

\def\ran{\rangle}

\def\TY(#1,#2,#3){#1^{(#2)}_{#3}}

\def\xxi(#1,#2,#3){\displaystyle {}^{#1}\Xi^{(#2)}_{#3}}
\def\xsi(#1,#2,#3){\displaystyle {}^{#1}\Sigma^{(#2)}_{#3}}
\def\xE(#1,#2,#3){\displaystyle {}^{#1}E_{#2}[#3]}
\def\xF(#1,#2){\displaystyle {}^{#1}F_{#2}}
\def\xx(#1,#2){\displaystyle {}^{#1}\Xi_{#2}}
\def\W1{W(\varpi_1)}

\def\m@th{\mathsurround=0pt}
\def\fsquare(#1,#2){
\hbox{\vrule$\hskip-0.4pt\vcenter to #1{\normalbaselines\m@th
\hrule\vfil\hbox to #1{\hfill$\scriptstyle #2$\hfill}\vfil\hrule}$\hskip-0.4pt
\vrule}}

\newcommand{\ba}{\begin{array}}
\newcommand{\ea}{\end{array}}

\newcommand{\eq}{\begin{eqnarray}}
\newcommand{\eneq}{\end{eqnarray}}

\title{\textbf{\large{Monomial realizations and LS paths of
fundamental representations for rank $2$ Kac-Moody algebras
}}}
\author{\normalsize{YUKI KANAKUBO\thanks{Faculty of Basic Natural Science, Ibaraki University : {j\_chi\_sen\_you\_ky@eagle.sophia.ac.jp
}}}
}
\date{}

\begin{document}

\maketitle
\vspace{-10pt}

\begin{abstract}
For a Kac-Moody algebra $\mathfrak{g}$ of rank $2$ and a fundamental weight $\lambda$,
we explicitly give an isomorphism between the set of 
Lakshmibai-Seshadri paths $\mathbb{B}(\lambda)$ and
monomial realization $\mathcal{M}(\lambda)$.  
As an application, we also give an explicit form of monomial realization $\mathcal{M}(\lambda)$ in terms of Weyl groups.
\end{abstract}

\section{Introduction}

For a symmetrizable Kac-Moody algebra $\mathfrak{g}$, the crystal bases powerful combinatorial tools to study
representations of quantum groups $U_q(\mathfrak{g})$ \cite{Kas90}. A bunch of
combinatorial descriptions of crystal bases have been developed.

In \cite{Nak03, Kas03}, monomial realizations of crystal bases $B(\lambda)$ of highest weight representations
are invented, which realize elements of $B(\lambda)$ in terms of Laurent monomials of double indexed variables $\{X_{s,i} | s\in\mathbb{Z},i\in I\}$. Here, $\lambda$ is a dominant integral weight and $I$ is the index set for simple roots of $\mathfrak{g}$. 
In particular, monomial realizations for rank $2$ Kac-Moody algebras are studied in \cite{JKKS}.
In \cite{Lit94, Lit95}, a piece-wise linear path $\pi:[0,1]\rightarrow P\otimes \mathbb{R}$ called
Lakshmibai-Seshadri path (LS path) are introduced. Here, $P$ is a weight lattice. 
The endpoint $\pi(1)$ is called a weight of $\pi$ and denoted by ${\rm wt}(\pi)$.
One can define
operators $f_k$ and $e_k$ on LS paths for $k\in I$ called root operators and
the set of LS paths of shape $\lambda$ denoted by $\mathbb{B}(\lambda)$ has a crystal structure by wt, $f_k$, $e_k$ and some
$\varepsilon$ and $\varphi$-functions.
It is shown in \cite{Kas96,Jos} that $\mathbb{B}(\lambda)$ is isomorphic to $B(\lambda)$.

In this paper, we will explicitly give an isomorphism between the set of LS paths $\mathbb{B}(\lambda)$
and monomial realization $\mathcal{M}(\lambda)$ when $\mathfrak{g}$
is of rank $2$ and $\lambda$ is a fundamental weight. Each LS path is defined by a sequence of rational numbers
$0=a_0<a_1<\cdots<a_r=1$ and Weyl group elements $\tau_1<\cdots<\tau_r$. The isomorphism is constructed
by these rational numbers and Weyl group elements (see Sect.\ref{main-sec}).
As an application, we also give an explicit form of monomial realization in terms of Weyl groups.

The paper is organized as follows.
In Sect.2, we review crystals. In Sect.3 and Sect.4, we review LS path models and monomial realizations.
We give our main results in Sect.5 and prove it in Sect.6.

\vspace{2mm}

\nd
{\bf Acknowledgements}
I would like to thank Daisuke Sagaki for his helpful comments
and fruitful discussions.

\section{Crystal}

Following \cite{Kas93}, let us review crystals.

\subsection{Notation}

Let $\mathfrak{g}$ be a symmetrizable Kac-Moody algebra with
Cartan subalgebra $\mathfrak{h}$, generalized Cartan matrix $(a_{i,j})_{i,j\in I}$ and index set $I=\{1,2,\cdots,
{\rm rank}\ \mathfrak{h}\}$.
Let $P$ and $P^{\vee}$ be the weight lattice and dual weight lattice with canonical pairing
$\lan,\ran : P^{\vee}\times P\rightarrow \mathbb{Z}$. Let $\{\alpha_i\}_{i\in I}\subset P$ and
$\{h_i\}_{i\in I}\subset P^{\vee}$ be the sets of simple roots and simple coroots. 
Let $P^+\subset P$ denote the set of dominant integral weights.
The fundamental weights
$\Lambda_i\in P^+$ ($\in I$) are defined as $\Lambda_i(h_j)=\delta_{i,j}$.
We define $W$ as its Weyl group with simple reflections $s_i$ ($i\in I$).
For $a,b\in\mathbb{Z}$ with $a\leq b$, let $[a,b]_{\rm int}:=\{a,a+1,\cdots,b\}$ and for $x,y\in\mathbb{R}$ with $x\leq y$,
let $[x,y]$ denote the closed interval from $x$ to $y$.

\subsection{Crystals}

\begin{defn}\label{pcry}\cite{Kas93}
A set $\mathcal{B}$ equipped with the maps
${\rm wt}:\mathcal{B}\rightarrow P$,
$\varepsilon_k,\varphi_k:\mathcal{B}\rightarrow \mathbb{Z}\sqcup \{-\infty\}$
and $\tilde{e}_k$,$\tilde{f}_k:\mathcal{B}\rightarrow \mathcal{B}\sqcup\{0\}$
($k\in I$)
is said to be a $P$-{\it crystal} if 
the following holds: For $b,b'\in\mathcal{B}$ and $k\in I$,
\begin{itemize}
\item $\varphi_k(b)=\varepsilon_k(b)+\langle h_k,{\rm wt}(b)\rangle$,
\item \[
{\rm wt}(\tilde{e}_kb)={\rm wt}(b)+\alpha_k \text{ if }\tilde{e}_k(b)\in\mathcal{B},
\quad {\rm wt}(\tilde{f}_kb)={\rm wt}(b)-\alpha_k\text{ if }\tilde{f}_k(b)\in\mathcal{B},
\]
\[
\varepsilon_k(\tilde{e}_k(b))=\varepsilon_k(b)-1,\ \ 
\varphi_k(\tilde{e}_k(b))=\varphi_k(b)+1\text{ if }\tilde{e}_k(b)\in\mathcal{B},
\] 
\[
\varepsilon_k(\tilde{f}_k(b))=\varepsilon_k(b)+1,\ \ 
\varphi_k(\tilde{f}_k(b))=\varphi_k(b)-1\text{ if }\tilde{f}_k(b)\in\mathcal{B},
\] 
\item $\tilde{f}_k(b)=b'$ is equivalent to $b=\tilde{e}_k(b')$,
\item if $\varphi_k(b)=-\infty$ then it holds $\tilde{e}_k(b)=\tilde{f}_k(b)=0$.
\end{itemize}
Here, $0$ and $-\infty$ are additional elements which do not belong to $\mathcal{B}$ and $\mathbb{Z}$, respectively.
We call the maps $\tilde{e}_k$,$\tilde{f}_k$ {\it Kashiwara operators}.
\end{defn}
The crystal bases $B(\lambda)$ $(\lambda\in P^+)$ and $B(\infty)$ have $P$-crystal structures.

\begin{defn}
We define $P_{\rm cl}:=P/\{\lambda\in P|\lambda(h_j)=0\text{ for all }j\in I\}$.
A set $\mathcal{B}$ equipped with the maps
${\rm wt}:\mathcal{B}\rightarrow P_{\rm cl}$,
$\varepsilon_k,\varphi_k:\mathcal{B}\rightarrow \mathbb{Z}\sqcup \{-\infty\}$
and $\tilde{e}_k$,$\tilde{f}_k:\mathcal{B}\rightarrow \mathcal{B}\sqcup\{0\}$
($k\in I$)
is said to be a $P_{\rm cl}$-{\it crystal} if the conditions in Definition \ref{pcry} hold.
\end{defn}

\begin{rem}
\begin{enumerate}
\item[(1)]
A $P$-crystal $\mathcal{B}$ naturally has a $P_{\rm cl}$-crystal structure by defining ${\rm wt}:\mathcal{B}\rightarrow P_{\rm cl}$
as the composition of ${\rm wt}:\mathcal{B}\rightarrow P$ and the projection $P\rightarrow P_{\rm cl}$.
\item[(2)]
By an abuse of notation, we use the same symbol $\Lambda_j$, $\alpha_j\in P_{\rm cl}$ for the image of
$\Lambda_j$, $\alpha_j\in P$ under the projection $P\rightarrow P_{\rm cl}$.
\end{enumerate}
\end{rem}

\begin{defn}
Let $\mathcal{B}_1$ and $\mathcal{B}_2$ be $P$-crystals (resp. $P_{\rm cl}$-crystals).
\begin{enumerate}
\item[$(1)$]
A map $\psi : \mathcal{B}_1\sqcup\{0\}\rightarrow \mathcal{B}_2\sqcup\{0\}$
is called a {\it strict morphism} 
if $\psi(0)=0$ and the following holds:
\begin{itemize}
\item For $k\in I$, if $b\in\mathcal{B}_1$ and $\psi(b)\in \mathcal{B}_2$ then we have
\[
{\rm wt}(\psi(b))={\rm wt}(b),\quad
\varepsilon_k(\psi(b))=\varepsilon_k(b),\quad
\varphi_k(\psi(b))=\varphi_k(b),
\]
\item $\tilde{e}_k(\psi(b))=\psi(\tilde{e}_k(b))$ and $\tilde{f}_k(\psi(b))=\psi(\tilde{f}_k(b))$ for $k\in I$ and $b\in\mathcal{B}_1$, where
we set $\tilde{e}_k(0)=\tilde{f}_k(0)=0$.
\end{itemize}
\item[$(2)$] If a strict morphism 
$\psi : \mathcal{B}_1\sqcup\{0\}\rightarrow \mathcal{B}_2\sqcup\{0\}$ is a bijection then
$\psi$ is said to be an {\it isomorphism}. Then we say $\mathcal{B}_1$ is isomorphic to $\mathcal{B}_2$ as $P$-crystals 
(resp. $P_{\rm cl}$-crystals).
\end{enumerate}
\end{defn}

\section{LS paths}

\subsection{Definition of LS paths}

Following \cite{Lit94}, let us review the Lakshmibai-Seshadri paths.
For $\lambda\in P^+$, let $W_{\lambda}:=\{w\in W|w\lambda=\lambda\}$ be the stabilizer of $\lambda$
and let $\geq$ be the Bruhat order on $W/W_{\lambda}$.

Let $\underline{\tau}=(\tau_1,\cdots,\tau_r)$ be a sequence of elements in $W/W_{\lambda}$ such that
$\tau_1>\cdots>\tau_r$. Let $\underline{a}=(a_0,a_1,\cdots,a_r)$ be a sequence of rational numbers such that
$a_0:=0<a_1<\cdots<a_r:=1$. The pair $\pi=(\underline{\tau};\underline{a})$ is called a rational $W$-path of shape $\lambda$.
The path $\pi$ is identified with a piecewise linear path
$\pi :[0,1]\rightarrow P_{\mathbb{R}}$ defined as
\[
\pi(t):=
\sum_{i=1}^{j-1}(a_i-a_{i-1})\tau_i(\lambda)
+(t-a_{j-1})\tau_j(\lambda)\quad \text{for }a_{j-1}\leq t\leq a_j,
\]
where $P_{\mathbb{R}}:=P\otimes_{\mathbb{Z}} \mathbb{R}$. Let $h_{k}:[0,1]\rightarrow\mathbb{R}$ be a function
defined as $t\mapsto\lan h_k,\pi(t)\ran$ for $k\in I$.

For $\tau,\sigma\in W/W_{\lambda}$ such that $\tau>\sigma$ and $a\in\mathbb{Q}$ such that $0<a<1$,
one defines an $a$-chain for the pair $(\tau,\sigma)$ as
\[
\kappa_0:=\tau>\kappa_1:=s_{\beta_1}\tau>\kappa_2:=s_{\beta_2}s_{\beta_1}\tau>\cdots>\kappa_s:=s_{\beta_s}\cdots s_{\beta_1}\tau=\sigma,
\]
where
$\beta_1,\cdots,\beta_s$ are positive real roots satisfying for each $i=1,\cdots,s$,
\[
l(\kappa_{i})=l(\kappa_{i-1})-1 \quad \text{and }a\lan \beta_i^{\vee},\kappa_i(\lambda)\ran \in\mathbb{Z}.
\]
\begin{defn}\cite{Lit94}
Let $\pi=(\underline{\tau};\underline{a})$ be
a rational $W$-path of shape $\lambda$ with
$\underline{\tau}=(\tau_1,\cdots,\tau_r)$ and $\underline{a}=(a_0,a_1,\cdots,a_r)$.
We say $\pi$ is a {\it Lakshmibai-Seshadri path} (in short LS path) of shape $\lambda$ if for each $i=1,2,\cdots,r-1$,
there is an $a_i$-chain for the pair $(\tau_i,\tau_{i+1})$. Let
$\mathbb{B}(\lambda)$ denote the set of Lakshmibai-Seshadri paths of shape $\lambda$. 
\end{defn}

\subsection{Action of Kashiwara operators on LS paths}\label{3-2}

We fix an index $k\in I$ and a LS path $\pi=(\underline{\tau};\underline{a})=(\tau_1,\cdots,\tau_r;a_0,a_1,\cdots,a_r)$. 
Let
\begin{equation}\label{Qdef}
Q:={\rm min}\{h_{k}(t) | t\in[0,1]\},\quad R:=h_{k}(1)-Q
\end{equation}
and $q$ (resp. $p$) be the minimal (resp. maximal) number in $\{0,1,2,\cdots,r\}$ such that $h_{k}(a_q)=Q$
(resp. $h_{k}(a_p)=Q$). It is shown $Q\in\mathbb{Z}$ \cite{Lit94}.
If $Q<0$ then we define
\[
y:={\rm max}\{s\in [0,q]_{\rm int}| h_{k}(t)\geq Q+1 \text{ for }t\in[0,a_{s}]\}.
\]
If $R>0$ then we also define
\[
x:={\rm min}\{s\in [p,r]_{\rm int}| h_{k}(t)\geq Q+1 \text{ for }t\in[a_{s},1]\}.
\]
\begin{defn}\label{opdef}\cite{Lit94}
\begin{enumerate}\item
If $R>0$ then $f_{k}(\pi)$ equals the following LS path:
\[
\begin{cases}
(\tau_1,\cdots,\tau_{p-1},s_k\tau_{p+1},\cdots,s_k\tau_x,\tau_{x+1},\cdots,\tau_r;a_0,\cdots,a_{p-1},a_{p+1},\cdots,a_r)  & \text{ if }\begin{array}{l}h_k(a_x)=Q+1,\\ s_{k}\tau_{p+1}=\tau_p,\end{array}\\
(\tau_1,\cdots,\tau_{p},s_k\tau_{p+1},\cdots,s_k\tau_x,\tau_{x+1},\cdots,\tau_r;a_0,\cdots,a_r)  & \text{ if }\begin{array}{l}h_k(a_x)=Q+1,\\ s_{k}\tau_{p+1}<\tau_p,\end{array}\\
(\tau_1,\cdots,\tau_{p-1},s_k\tau_{p+1},\cdots,s_k\tau_x,\tau_{x},\cdots,\tau_r;a_0,\cdots,a_{p-1},a_{p+1},\cdots,a_{x-1},a,a_{x},\cdots,a_r)  & \text{ if }\begin{array}{l}h_k(a_x)>Q+1,\\ s_{k}\tau_{p+1}=\tau_p,\end{array}\\
(\tau_1,\cdots,\tau_{p},s_k\tau_{p+1},\cdots,s_k\tau_x,\tau_{x},\cdots,\tau_r;a_0,\cdots,a_{x-1},a,a_{x},\cdots,a_r)  & \text{ if }\begin{array}{l}h_k(a_x)>Q+1,\\ s_{k}\tau_{p+1}<\tau_p.\end{array}
\end{cases}
\]
Here, $a$ satisfies $a_{x-1}<a<a_x$ and $h_k(a)=Q+1$. When $p=0$, we understand $s_k\tau_{p+1}<\tau_p$.
\item If $Q<0$ then $e_{k}(\pi)$ equals the following LS path:
\[
\begin{cases}
(\tau_1,\cdots,\tau_{y},s_k\tau_{y+1},\cdots,s_k\tau_q,\tau_{q+2},\cdots,\tau_r;a_0,\cdots,a_{q},a_{q+2},\cdots,a_r)  & \text{ if }\begin{array}{l}h_k(a_y)=Q+1,\\ s_{k}\tau_{q}=\tau_{q+1},\end{array}\\
(\tau_1,\cdots,\tau_{y},s_k\tau_{y+1},\cdots,s_k\tau_q,\tau_{q+1},\cdots,\tau_r;a_0,\cdots,a_r)  & \text{ if }\begin{array}{l}h_k(a_y)=Q+1,\\ s_{k}\tau_{q}>\tau_{q+1},\end{array}\\
(\tau_1,\cdots,\tau_{y+1},s_k\tau_{y+1},\cdots,s_k\tau_q,\tau_{q+2},\cdots,\tau_r;a_0,\cdots,a_{y},a,a_{y+1},\cdots,a_{q},a_{q+2},\cdots,a_r)  & \text{ if }\begin{array}{l}h_k(a_y)>Q+1,\\ s_{k}\tau_{q}=\tau_{q+1},\end{array}\\
(\tau_1,\cdots,\tau_{y+1},s_k\tau_{y+1},\cdots,s_k\tau_q,\tau_{q+1},\cdots,\tau_r;a_0,\cdots,a_{y},a,a_{y+1},\cdots,a_r)  & \text{ if }\begin{array}{l}h_k(a_y)>Q+1,\\ s_{k}\tau_{q}>\tau_{q+1}.\end{array}\\
\end{cases}
\]
Here, $a$ satisfies $a_{y}<a<a_{y+1}$ and $h_k(a)=Q+1$. When $q=r$, we understand $s_k\tau_q>\tau_{q+1}$.
\end{enumerate}

\end{defn}

We define
\begin{equation}\label{LSwt}
{\rm wt}(\pi):=\pi(1)=\sum_{i=1}^{r}(a_i-a_{i-1})\tau_i(\lambda),
\end{equation}
\begin{equation}\label{LSep}
\varepsilon_k(\pi):=\mathbb{Z}\cap
{\rm max}\{-h_{k}(t) | t\in[0,1]\},\quad
\varphi_k(\pi):=\lan h_k,\pi(1)\ran+\varepsilon_k(\pi).
\end{equation}
For $\pi\in\mathbb{B}(\lambda)$, it holds
$\varepsilon_k(\pi)=-Q$ and $\varphi_k(\pi)=R$
with $Q$, $R$ in (\ref{Qdef}).

\begin{thm}\cite{Kas96,Jos}
For a dominant integral weight $\lambda$,
the set $\mathbb{B}(\lambda)$ with the above maps is a $P$-crystal isomorphic
to $B(\lambda)$.
\end{thm}

\subsection{Rank $2$ case}\label{rank2}

We assume $\mathfrak{g}$ is rank $2$. The generalized Cartan matrix $A$ is in the form
\[
A=
\begin{pmatrix}
2 & -a \\
-b & 2
\end{pmatrix}
\]
with $a,b\in\mathbb{Z}_{\geq0}$.

One defines
\[
\tau_{2n}^+=(s_2s_1)^n,\quad
\tau_{2n+1}^+=s_1(s_2s_1)^n,\quad
\tau_{2n}^-=(s_1s_2)^n,\quad
\tau_{2n+1}^-=s_2(s_1s_2)^n
\]
for $n\in \mathbb{Z}_{\geq0}$. We see that
\[
W/W_{\Lambda_i}
=\{\tau_k^{\varepsilon(i)} | 0\leq k<N\},
\]
where $N:=|\lan s_2s_1\ran|$ and 
\[
\varepsilon(i)=
\begin{cases}
+ & \text{ if }i=1,\\
- & \text{ if }i=2.
\end{cases}
\]
Note that $N=\infty$ when $\mathfrak{g}$ is not finite dimensional.
Let $\xi$ be one of solutions for $x^2-(ab-2)x+1=0$ and
\[
[n]:=\xi^n+\xi^{n-2}+\cdots+\cdots+\xi^{-n+2}+\xi^{-n}
\]
for $n\geq0$ and $[-1]:=0$.
We also set
\[
d_{2n}^+:=-b[n-1],\quad
d_{2n+1}^+:=-[n]-[n-1],\quad
d_{2n}^-:=-a[n-1],\quad
d_{2n+1}^-:=-[n]-[n-1].
\]

\begin{thm}\label{Sagthm}\cite{Sag}
\[
\mathbb{B}(\Lambda_i)
=
\left\{
(\tau_{m_2}^{\varepsilon(i)},\tau_{m_2-1}^{\varepsilon(i)},\cdots,\tau_{m_1+1}^{\varepsilon(i)},\tau_{m_1}^{\varepsilon(i)};
0,a_{m_2},a_{m_2-1},\cdots,a_{m_1+1},1)|
\begin{array}{l}
N> m_2\geq m_1\geq0,\\
0<a_{m_2}<\cdots<a_{m_1+1}<1,\\
a_jd_j^{\varepsilon(i)}\in\mathbb{Z}\text{ for }j=m_1+1,\cdots,m_2
\end{array}
\right\}.
\]
\end{thm}

\begin{ex}\label{LS-B2-ex1}
We consider the case the Cartan matrix is of type $B_2$:
$A=
\begin{pmatrix}
2 & -1 \\
-2 & 2
\end{pmatrix}$. 
One obtains
\[
W/W_{\Lambda_1}=\{e,s_1,s_2s_1,s_1s_2s_1\}.
\]
Using the notation as above, it holds $\{\xi,\xi^{-1}\}=\{\sqrt{-1},-\sqrt{-1}\}$ so that
\[
[n]=
\begin{cases}
0 & \text{ if }n\text{ is odd},\\
-1 & \text{ if }n\equiv 2\ \text{mod }4,\\
1 & \text{ if }n\equiv 0\ \text{mod }4.
\end{cases}
\]
Hence, $d_1^+=-1$, $d_2^+=-2$, $d_3^+=-1$ so that
$a_1=a_3=1$ and $a_2\in\{1,\frac{1}{2}\}$ if $a_jd_j^+\in\mathbb{Z}$ for $j=1,2,3$. Therefore, when $m_2>m_1$, if $m_1=0$ (resp. $m_2\geq3$)
then it contradicts $0<a_1<1$ (resp. $0<a_3<1$).
Thus,
if 
$
(\tau_{m_2}^{+},\tau_{m_2-1}^{+},\cdots,\tau_{m_1+1}^{+},\tau_{m_1}^{+};
0,a_{m_2},a_{m_2-1},\cdots,a_{m_1+1},1)\in
\mathbb{B}(\Lambda_1)$ then it holds either $0\leq m_1=m_2\leq 3$ or $m_2=2$, $m_1=1$ and $a_2=\frac{1}{2}$, that is,
\[
\mathbb{B}(\Lambda_1)
=
\{
(e;0,1),\ (s_1;0,1),\ (s_2s_1;0,1),\ (s_1s_2s_1;0,1),\ (s_2s_1,s_1;0,\frac{1}{2},1) 
\}.
\]
The crystal graph is as follows:
\[
(e;0,1)\overset{1}{\rightarrow}(s_1;0,1)
\overset{2}{\rightarrow} (s_2s_1,s_1;0,\frac{1}{2},1) \overset{2}{\rightarrow}(s_2s_1;0,1)
\overset{1}{\rightarrow} (s_1s_2s_1;0,1).
\]
\end{ex}

\begin{ex}\label{LS-A11-ex1}
We consider the case the generalized Cartan matrix is of type $A^{(1)}_1$:
$A=
\begin{pmatrix}
2 & -2 \\
-2 & 2
\end{pmatrix}$. 
A part of the crystal graph for $\mathbb{B}(\Lambda_1)$ is as follows:
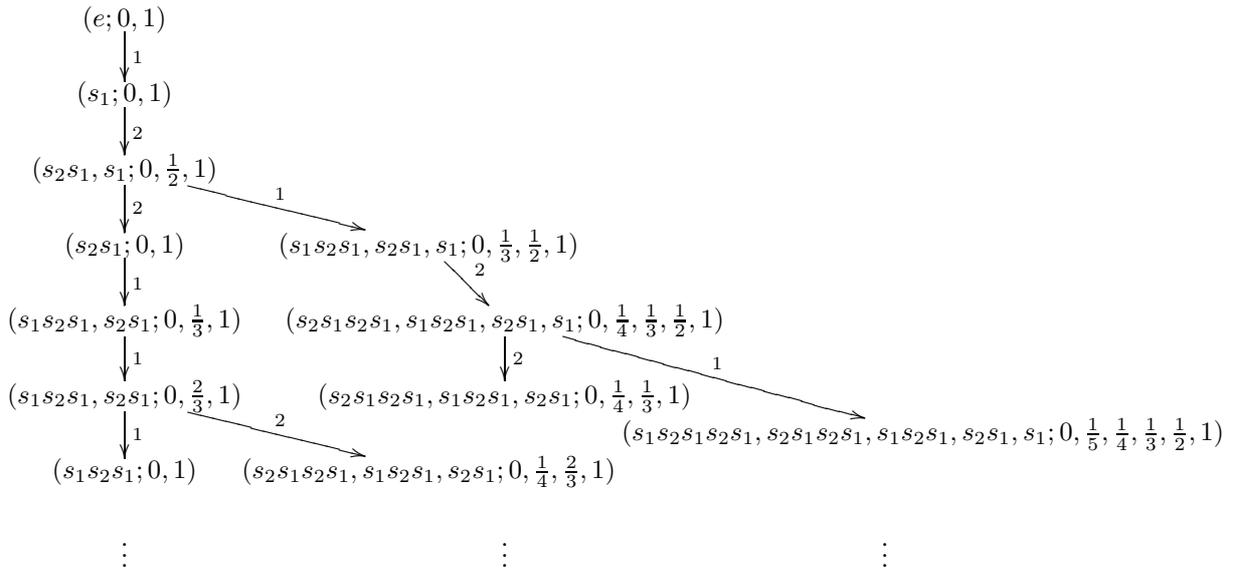
\[
\begin{xy}
(0,0) *{(e;0,1)}="1",
(0,-10) *{(s_1;0,1)}="2",
(0,-20) *{(s_2s_1,s_1;0,\frac{1}{2},1)}="3",
(0,-30) *{(s_2s_1;0,1)}="4-1",
(40,-30) *{(s_1s_2s_1,s_2s_1,s_1;0,\frac{1}{3},\frac{1}{2},1)}="4-2",
(0,-40) *{(s_1s_2s_1,s_2s_1;0,\frac{1}{3},1)}="5-1",
(50,-40) *{(s_2s_1s_2s_1,s_1s_2s_1,s_2s_1,s_1;0,\frac{1}{4},\frac{1}{3},\frac{1}{2},1)}="5-2",
(0,-50) *{(s_1s_2s_1,s_2s_1;0,\frac{2}{3},1)}="6-1",
(50,-50) *{(s_2s_1s_2s_1,s_1s_2s_1,s_2s_1;0,\frac{1}{4},\frac{1}{3},1)}="6-2",
(105,-55) *{(s_1s_2s_1s_2s_1,s_2s_1s_2s_1,s_1s_2s_1,s_2s_1,s_1;0,\frac{1}{5},\frac{1}{4},\frac{1}{3},\frac{1}{2},1)}="6-3",
(0,-60) *{(s_1s_2s_1;0,1)}="7-1",
(40,-60) *{(s_2s_1s_2s_1,s_1s_2s_1,s_2s_1;0,\frac{1}{4},\frac{2}{3},1)}="7-2",
(0,-70) *{\vdots}="dot1",
(50,-70) *{\vdots}="dot2",
(100,-70) *{\vdots}="dot3",
\ar@{->} "1";"2"^{1}
\ar@{->} "2";"3"^{2}
\ar@{->} "3";"4-1"^{2}
\ar@{->} "3";"4-2"^{1}
\ar@{->} "4-1";"5-1"^{1}
\ar@{->} "4-2";"5-2"^{2}
\ar@{->} "5-1";"6-1"^{1}
\ar@{->} "5-2";"6-2"^{2}
\ar@{->} "5-2";"6-3"^{1}
\ar@{->} "6-1";"7-1"^{1}
\ar@{->} "6-1";"7-2"^{2}
\end{xy}
\]

\end{ex}

\section{Monomial Realizations}\label{monosec}

Following \cite{Kas03}, we recall the notion of monomial realizations for crystal bases of highest weight representations.
One defines a crystal structure on the set of Laurent monomials
\begin{gather*}
{\mathcal Y}:=\left\{X=\prod\limits_{s \in \mathbb{Z},\ i \in I}
X_{s,i}^{\zeta_{s,i}}\, \Bigg| \,\zeta_{s,i} \in \mathbb{Z},\ 
\zeta_{s,i}= 0\text{ except for finitely many }(s,i) \right\}.
\end{gather*}
For $X=\prod\limits_{s \in \mathbb{Z},\; i \in I} X_{s,i}^{\zeta_{s,i}}\in {\mathcal Y}$, we define
$\text{wt}(X):= \sum\limits_{i,s}\zeta_{s,i}\Lambda_i\in P_{\rm cl}$ and

\vspace{-3mm}

\begin{gather*}
\varphi_i(X):=\max\left\{\! \sum\limits_{k\leq s}\zeta_{k,i}  \,|\, s\in \mathbb{Z} \!\right\},\!
\quad
\varepsilon_i(X):=\varphi_i(X)-\text{wt}(X)(h_i).
\end{gather*}
Next, we put
\[
A_{s,k}:=X_{s,k}X_{s+1,k}\prod\limits_{j\neq k,\ a_{j,k}<0}X_{s+p_{j,k},j}^{a_{j,k}},\ \ \ (s\in\mathbb{Z},\ k\in I)
\]
where $(a_{j,k})$ is a generalized Cartan matrix of $\mathfrak{g}$ and the integers $\{p_{j,k}\}_{j,k\in I; a_{j,k<0}}$ are defined as
\[
p_{j,k}=
\begin{cases}
1 & \text{if }k>j,\\
0 & \text{if }j>k.
\end{cases}
\]
Actions of Kashiwara operators are defined as follows:
\begin{gather*}
\tilde{f}_iX:=
\begin{cases}
A_{n_{f_i},i}^{-1}X & \text{if} \quad \varphi_i(X)>0,
\\
0 & \text{if} \quad \varphi_i(X)=0,
\end{cases}
\quad
\tilde{e}_iX:=
\begin{cases}
A_{n_{e_i},i}X & \text{if} \quad \varepsilon_i(X)>0,
\\
0 & \text{if} \quad \varepsilon_i(X)=0,
\end{cases}
\end{gather*}
where
\begin{gather*}
n_{f_i}:=\min \left\{r\in\mathbb{Z} \,\Bigg|\, \varphi_i(X)= \sum\limits_{k\leq r}\zeta_{k,i}\right\},
\qquad
n_{e_i}:=\max \left\{r\in\mathbb{Z} \,\Bigg|\, \varphi_i(X)= \sum\limits_{k\leq r}\zeta_{k,i}\right\}.
\end{gather*}
Then we get the following:
\begin{thm}\label{mono-real}\cite{Kas03,Nak03}
\begin{enumerate}
\item[(i)] The set ${\mathcal Y}$ with the above maps
${\rm wt}$, $\varepsilon_i$, $\varphi_i$ and $\tilde{e}_i$, $\tilde{f}_i$ $(i\in I)$ is a $P_{\rm cl}$-crystal.
\item[(ii)] For a monomial $X \in {\mathcal Y}$ such that $\tilde{\varepsilon}_i(X)=0$ for all $i \in I$, the set
\[
\{\tilde{f}_{j_m}\cdots\tilde{f}_{j_1}X | m\in\mathbb{Z}_{\geq0},\ j_1,\cdots,j_m \in I \}\setminus\{0\}
\]
is isomorphic to $B({\rm wt}(X))$ as $P_{\rm cl}$-crystals.
\end{enumerate}
\end{thm}

\begin{ex}\label{LS-B2-ex2}
We consider the case the Cartan matrix is of type $B_2$:
$A=
\begin{pmatrix}
2 & -1 \\
-2 & 2
\end{pmatrix}$. 

Then the crystal graph with the highest weight vector $X_{s,1}$, which is isomorphic to $B(\Lambda_1)$ is as follows:
\[
X_{s,1}\overset{1}{\rightarrow}\frac{X_{s,2}^2}{X_{s+1,1}}
\overset{2}{\rightarrow} \frac{X_{s,2}}{X_{s+1,2}} \overset{2}{\rightarrow} \frac{X_{s+1,1}}{X_{s+1,2}^2}
\overset{1}{\rightarrow} \frac{1}{X_{s+2,1}}.
\]
\end{ex}

\begin{ex}\label{LS-A11-ex2}
We consider the case the generalized Cartan matrix is of type $A^{(1)}_1$:
$A=
\begin{pmatrix}
2 & -2 \\
-2 & 2
\end{pmatrix}$.
Then the crystal graph with the highest weight vector $X_{s,1}$ is as follows:
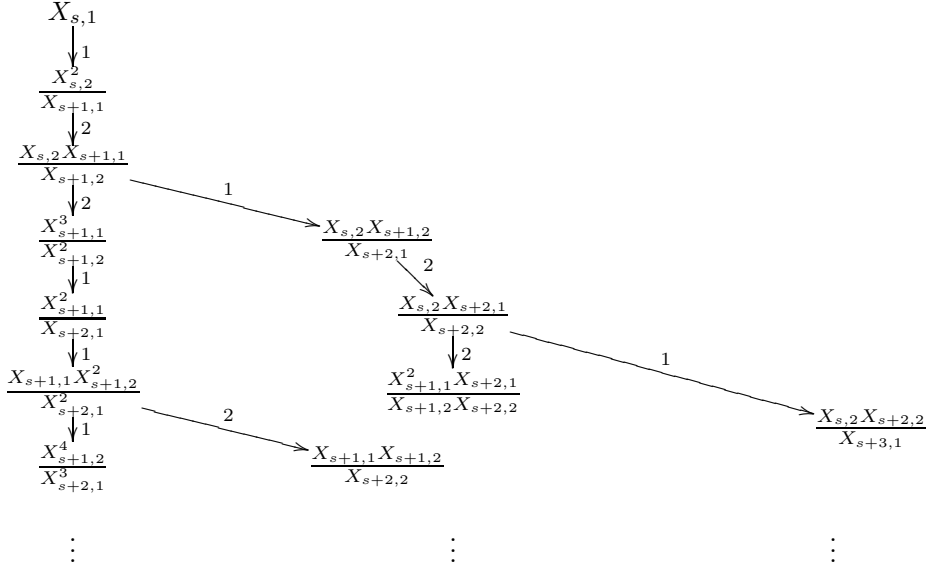
\[
\begin{xy}
(0,0) *{X_{s,1}}="1",
(0,-10) *{\frac{X_{s,2}^2}{X_{s+1,1}}}="2",
(0,-20) *{\frac{X_{s,2}X_{s+1,1}}{X_{s+1,2}}}="3",
(0,-30) *{\frac{X_{s+1,1}^3}{X_{s+1,2}^2}}="4-1",
(40,-30) *{\frac{X_{s,2}X_{s+1,2}}{X_{s+2,1}}}="4-2",
(0,-40) *{\frac{X_{s+1,1}^2}{X_{s+2,1}}}="5-1",
(50,-40) *{\frac{X_{s,2}X_{s+2,1}}{X_{s+2,2}}}="5-2",
(0,-50) *{\frac{X_{s+1,1}X_{s+1,2}^2}{X_{s+2,1}^2}}="6-1",
(50,-50) *{\frac{X_{s+1,1}^2X_{s+2,1}}{X_{s+1,2}X_{s+2,2}}}="6-2",
(105,-55) *{\frac{X_{s,2}X_{s+2,2}}{X_{s+3,1}}}="6-3",
(0,-60) *{\frac{X_{s+1,2}^4}{X_{s+2,1}^3}}="7-1",
(40,-60) *{\frac{X_{s+1,1}X_{s+1,2}}{X_{s+2,2}}}="7-2",
(0,-70) *{\vdots}="dot1",
(50,-70) *{\vdots}="dot2",
(100,-70) *{\vdots}="dot3",
\ar@{->} "1";"2"^{1}
\ar@{->} "2";"3"^{2}
\ar@{->} "3";"4-1"^{2}
\ar@{->} "3";"4-2"^{1}
\ar@{->} "4-1";"5-1"^{1}
\ar@{->} "4-2";"5-2"^{2}
\ar@{->} "5-1";"6-1"^{1}
\ar@{->} "5-2";"6-2"^{2}
\ar@{->} "5-2";"6-3"^{1}
\ar@{->} "6-1";"7-1"^{1}
\ar@{->} "6-1";"7-2"^{2}
\end{xy}
\]

\end{ex}

\section{Main results}\label{main-sec}

Let $\mathfrak{g}$ be a rank $2$ Kac-Moody algebra. 
We fix $s\in\mathbb{Z}$ and $i\in I$ and define $i'\in I$ as $\{i,i'\}=\{1,2\}=I$.
Let $\mathcal{M}(\Lambda_i)$ be the monomial realization of $B(\Lambda_i)$ with the highest weight vector
$X_{s,i}$ in Theorem \ref{mono-real}, that is,
\[
\mathcal{M}(\Lambda_i)=
\{\tilde{f}_{j_m}\cdots\tilde{f}_{j_1}X_{s,i} | m\in\mathbb{Z}_{\geq0},\ j_1,\cdots,j_m \in I \}\setminus\{0\},
\]
which is isomorphic to $B(\Lambda_i)$. For $w\in W$, let $v_{w\Lambda_i}$ denote the extremal weight vector in $B(\Lambda_i)$ with
weight $w\Lambda_{i}\in P$.

For $m\in\mathbb{Z}_{\geq0}$
such that $2m\leq N-1$,
we define integers $P_{m,i}$, $P_{m,i'}$ as follows:
\[
(s_{i'}s_i)^m\Lambda_i=P_{m,i}\Lambda_i-P_{m,i'}\Lambda_{i'} \text{ in }P_{\rm cl}.
\]
Note that when $2m+1\leq N-1$, it holds
\[
s_i(s_{i'}s_i)^m\Lambda_i=P_{m+1,i'}\Lambda_{i'}-P_{m,i}\Lambda_{i},
\]
where if $2m+2\geq N$ then we define $P_{m+1,i'}=0$.
For $w=(s_{i'}s_i)^m\in W$ as above, we define the following elements in $\mathcal{Y}$:
\begin{equation}\label{pr-0}
M(v_{w\Lambda_i}):=M_s(v_{w\Lambda_i})=
\begin{cases}
\left(
\frac{X^{P_{m,1}}_{s+m,1}}{X^{P_{m,2}}_{s+m,2}}
\right) & {\rm if}\ i=1,\\
\left(
\frac{X^{P_{m,2}}_{s+m,2}}{X^{P_{m,1}}_{s+m+1,1}}
\right) & {\rm if}\ i=2,
\end{cases}
\end{equation}
\begin{equation}\label{pr-00}
M(v_{s_iw\Lambda_i}):=M_s(v_{s_iw\Lambda_i})=
\begin{cases}
\left(
\frac{X^{P_{m+1,2}}_{s+m,2}}{X^{P_{m,1}}_{s+m+1,1}}
\right) & {\rm if}\ i=1,\\
\left(
\frac{X^{P_{m+1,1}}_{s+m+1,1}}{X^{P_{m,2}}_{s+m+1,2}}
\right) & {\rm if}\ i=2.
\end{cases}
\end{equation}
These are extremal weight vectors in $\mathcal{M}(\Lambda_i)$ with weight $w\Lambda_i$ and $s_iw\Lambda_i$, respectively.

We define a map
\[
\Phi_s=\Phi : \mathbb{B}(\Lambda_i) \rightarrow \mathcal{Y}
\]
as
\[
(\tau_1,\cdots,\tau_r ; a_0=0,a_1,\cdots,a_{r-1},a_r=1)
\mapsto
\prod^{r}_{l=1}
\left(M_s(v_{\tau_l\Lambda_i}) \right)^{a_l-a_{l-1}}.
\]
\begin{thm}\label{mainthm}
The above map is well-defined and  induces an isomorphism between $\mathbb{B}(\Lambda_i)$ and $\mathcal{M}(\Lambda_i)$
as $P_{\rm cl}$-crystals.
\end{thm}
Here, well-defined means the image of $\Phi$ is included in $\mathcal{Y}$.
Combining with Theorem \ref{Sagthm}, we get the following expression of Monomial realization $\mathcal{M}(\Lambda_i)$.
\begin{cor}
In the notation of subsection \ref{rank2}, we have
\[
\mathcal{M}(\Lambda_i)
=\left\{\prod^{m_2}_{l=m_1}
\left(M(v_{\tau_l^{\varepsilon(i)}\Lambda_i}) \right)^{a_l-a_{l+1}} | 
\begin{array}{l}
N> m_2\geq m_1\geq0,\\
a_{m_2+1}:=0<a_{m_2}<\cdots<a_{m_1+1}<a_{m_1}:=1,\\
a_jd_j^{\varepsilon(i)}\in\mathbb{Z}\text{ for }j=m_1+1,\cdots,m_2
\end{array}
\right\}.
\]
\end{cor}

\begin{ex}
We consider the case the Cartan matrix is of type $B_2$:
$A=
\begin{pmatrix}
2 & -1 \\
-2 & 2
\end{pmatrix}$. It is the same setting as in Example \ref{LS-B2-ex1} and \ref{LS-B2-ex2}.
We obtained
\[
\mathbb{B}(\Lambda_1)
=
\{
(e;0,1),\ (s_1;0,1),\ (s_2s_1;0,1),\ (s_1s_2s_1;0,1),\ (s_2s_1,s_1;0,\frac{1}{2},1) 
\}.
\] 
By the isomorphism $\Phi$, each element goes to an element of $\mathcal{M}(\Lambda_1)$:
\[
(e;0,1)\mapsto X_{s,1},\quad (s_1;0,1)\mapsto \frac{X_{s,2}^2}{X_{s+1,1}},\quad
(s_2s_1;0,1)\mapsto \frac{X_{s+1,1}}{X_{s+1,2}^2},\quad 
(s_1s_2s_1;0,1)\mapsto \frac{1}{X_{s+2,1}},
\]
\[
(s_2s_1,s_1;0,\frac{1}{2},1) 
\mapsto
M(v_{s_2s_1\Lambda_1})^{\frac{1}{2}}
M(v_{s_1\Lambda_1})^{\frac{1}{2}}
=
\left(
\frac{X_{s,2}^2}{X_{s+1,1}}
\right)^{\frac{1}{2}}
\left(
\frac{X_{s+1,1}}{X_{s+1,2}^2}
\right)^{\frac{1}{2}}
=\frac{X_{s,2}}{X_{s+1,2}}.
\]
Hence, the monomials coincide with those in
Example \ref{LS-B2-ex2}.
\end{ex}

\begin{ex}
We consider the case the generalized Cartan matrix is of type $A^{(1)}_1$:
$A=
\begin{pmatrix}
2 & -2 \\
-2 & 2
\end{pmatrix}$, which is the same setting as in Example \ref{LS-A11-ex1} and \ref{LS-A11-ex2}.
By the isomorphism $\Phi$, the paths in Example \ref{LS-A11-ex1} go to the corresponding monomials
in Example \ref{LS-A11-ex2}.
For example, the element 
\[
(s_1s_2s_1s_2s_1,s_2s_1s_2s_1,s_1s_2s_1,s_2s_1,s_1;0,\frac{1}{5},\frac{1}{4},\frac{1}{3},\frac{1}{2},1)
\]
goes to
\begin{eqnarray*}
& &M(v_{s_1s_2s_1s_2s_1\Lambda_1})^{\frac{1}{5}-0}M(v_{s_2s_1s_2s_1\Lambda_1})^{\frac{1}{4}-\frac{1}{5}}
M(v_{s_1s_2s_1\Lambda_1})^{\frac{1}{3}-\frac{1}{4}}
M(v_{s_2s_1\Lambda_1})^{\frac{1}{2}-\frac{1}{3}}
M(v_{s_1\Lambda_1})^{1-\frac{1}{2}}\\
&=&\left(\frac{X_{s+2,2}^6}{X_{s+3,1}^5}\right)^{\frac{1}{5}}
\left(\frac{X_{s+2,1}^5}{X_{s+2,2}^4}\right)^{\frac{1}{20}}
\left(\frac{X_{s+1,2}^4}{X_{s+2,1}^3}\right)^{\frac{1}{12}}
\left(\frac{X_{s+1,1}^3}{X_{s+1,2}^2}\right)^{\frac{1}{6}}
\left(\frac{X_{s,2}^2}{X_{s+1,1}}\right)^{\frac{1}{2}}\\
&=&\frac{X_{s,2}X_{s+2,2}}{X_{s+3,1}}.
\end{eqnarray*}
\end{ex}

\begin{rem}
We can not apply our result to Kac-Moody algebras of other ranks and
other weights.  
\end{rem}

\section{Proof}

[Proof of Theorem \ref{mainthm}]

\vspace{1mm}

\underline{Well-definedness and $\varepsilon_k(\Phi(\pi))=\varepsilon_k(\pi)$}

\vspace{2mm}

We fix $k\in\{1,2\}$.
For $\pi=(\tau_1,\cdots,\tau_r ; a_0=0,a_1,\cdots,a_{r-1},a_r=1)\in\mathbb{B}(\Lambda_i)$
and $j\in[0,r]_{\rm int}$, we put
\[
y_j:=\sum^{j}_{l=1} \langle h_k,\tau_l\Lambda_i\rangle (a_l-a_{l-1})=h_k(a_j).
\]
In particular, it holds $y_r=\langle h_k,{\rm wt}(\pi)\rangle\in\mathbb{Z}$ by (\ref{LSwt}). Let us prove the exponent of variables
in the form $X_{t,k}$ ($t\in\mathbb{Z}$) in $\Phi(\pi)$ is an integer and $\varepsilon_k(\Phi(\pi))=\varepsilon_k(\pi)$.

\vspace{2mm}

Since $l(\tau_j)=l(\tau_{j+1})+1$ by Theorem \ref{Sagthm} and there exists an $a_j$-chain for $(\tau_j,\tau_{j+1})$ for $j\in[1,r-1]_{\rm int}$, we see that for $K\in\{1,2\}$,
\begin{equation}\label{pr-1}
a_j\langle h_K,\tau_{j}\Lambda_i \rangle
=-a_j\langle h_K,\tau_{j+1}\Lambda_i \rangle
\in\mathbb{Z}\ \ {\rm if}\ s_K\tau_j=\tau_{j+1}.
\end{equation}
If $j<r$ and $s_k\tau_j=\tau_{j+1}$ then
\begin{eqnarray*}
y_j&=&\sum^{j}_{l=1} \langle h_k,\tau_l\Lambda_i\rangle (a_l-a_{l-1})\\
&=& \langle h_k,\tau_j\Lambda_i\rangle a_j 
+ \sum^{j-1}_{l=1}\langle h_k,\tau_l\Lambda_i-\tau_{l+1}\Lambda_i\rangle a_l\\
&=& \langle h_k,\tau_j\Lambda_i\rangle a_j 
+ \sum_{l\in[1,j-1]_{\rm int},\ s_k\tau_{l}=\tau_{l+1}}\langle h_k,\tau_l\Lambda_i-\tau_{l+1}\Lambda_i\rangle a_l\\
& &+ \sum_{l\in[1,j-1]_{\rm int},\ s_{k'}\tau_{l}=\tau_{l+1}}\langle h_k,\tau_l\Lambda_i-\tau_{l+1}\Lambda_i\rangle a_l,
\end{eqnarray*}
where $k'\in\{1,2\}\setminus \{k\}$. By (\ref{pr-1}), the first term and second term are in $\mathbb{Z}$.
As for the third term, if $s_{k'}\tau_{l}=\tau_{l+1}$ then
\begin{eqnarray*}
& &\langle h_k,\tau_l\Lambda_i-\tau_{l+1}\Lambda_i\rangle a_l\\
&=&\langle h_k,\tau_l\Lambda_i\rangle a_l - \langle h_k,\tau_{l+1}\Lambda_i\rangle a_l \\
&=&\langle h_k,\tau_l\Lambda_i\rangle a_l - \langle s_{k'}h_k,\tau_{l}\Lambda_i \rangle a_l \\
&=&\langle h_k,\tau_l\Lambda_i\rangle a_l - \langle h_k,\tau_{l}\Lambda_i\rangle a_l
+ \langle \alpha_{k'}(h_k)h_{k'},\tau_{l}\Lambda_i\rangle a_l\\
&=&\alpha_{k'}(h_k)\cdot a_l\cdot \langle h_{k'},\tau_{l}\Lambda_i \rangle.
\end{eqnarray*}
By (\ref{pr-1}), it holds $a_l\cdot \langle h_{k'},\tau_{l}\Lambda_i\rangle\in\mathbb{Z}$ so that 
$\langle h_k,\tau_l\Lambda_i-\tau_{l+1}\Lambda_i\rangle a_l\in\mathbb{Z}$. In this way,
we have 
\begin{equation}\label{pr-9}
y_j\in\mathbb{Z} \text{ if }s_k\tau_j=\tau_{j+1}\ \text{or }j=r.
\end{equation}

For $j\in[2,r]_{\rm int}$, 
we see that $-\langle h_k, \tau_j\Lambda_i \rangle\geq0$ if and only if
$\langle h_k, \tau_{j-1}\Lambda_i \rangle\geq0$. Thus, $M(v_{\tau_j\Lambda_i})$ has
a denominator $X_{t,k}^{-\langle h_k, \tau_j\Lambda_i \rangle}$ with
$-\langle h_k, \tau_j\Lambda_i \rangle\geq0$
in the expressions
(\ref{pr-0}) and (\ref{pr-00}) if and only if
$M(v_{\tau_{j-1}\Lambda_i})$ has a numerator 
$X_{t,k}^{\langle h_k, \tau_{j-1}\Lambda_i \rangle}$ with $\langle h_k, \tau_{j-1}\Lambda_i \rangle\geq0$.
Hence, the exponent of $X_{t,k}$ in $\Phi(\pi)$ is equal to 
\[
(a_{j-1}-a_{j-2})\langle h_k, \tau_{j-1}\Lambda_i \rangle+(a_j-a_{j-1})\langle h_k, \tau_j\Lambda_i \rangle
=y_{j}-y_{j-2},
\]
where $y_0=0$. In this case,
it holds that $s_k\tau_j=\tau_{j+1}$ (or $j=r$) and $s_k\tau_{j-2}=\tau_{j-1}$ if $j-2\geq1$ by Theorem \ref{Sagthm}. 
Thus, $y_j\in\mathbb{Z}$ and $y_{j-2}\in\mathbb{Z}$ by (\ref{pr-9}).
Therefore, the exponent of $X_{t,k}$ in $\Phi(\pi)$ is an integer.


If $M(v_{\tau_1\Lambda_i})$ has a denominator $X_{t,k}^{-\langle h_k, \tau_1\Lambda_i \rangle}$ with
$-\langle h_k, \tau_1\Lambda_i \rangle\geq0$
in the expressions
(\ref{pr-0}) and (\ref{pr-00}) then 
the exponent of $X_{t,k}$ in $\Phi(\pi)$ is equal to $(a_1-a_{0})\langle h_k, \tau_1\Lambda_i \rangle=y_1$.
In this case, it holds $s_k\tau_1=\tau_{2}$ or $r=1$ so that $y_1\in\mathbb{Z}$ by (\ref{pr-9}).

If $M(v_{\tau_r\Lambda_i})$ has a numerator 
$X_{t,k}^{\langle h_k, \tau_{r}\Lambda_i \rangle}$ with $\langle h_k, \tau_{r}\Lambda_i \rangle\geq0$
then the exponent of $X_{t,k}$ in $\Phi(\pi)$ is equal to
\[
(1-a_{r-1})\langle h_k, \tau_r\Lambda_i \rangle
=\langle h_k, \tau_r\Lambda_i \rangle+a_{r-1}\langle h_k, \tau_{r-1}\Lambda_i \rangle\in\mathbb{Z}.
\]
Therefore, the exponent of each variable $X_{t,k}$ in $\Phi(\pi)$ is an integer and $\Phi$ is well-defined.

It holds $\tau_{j+1}\Lambda_i(h_k)>0$
if and only if $y_j<y_{j+1}$ for $j\in[0,r-1]_{\rm int}$. Thus, it follows by (\ref{LSep}) that
\begin{eqnarray*}
\varepsilon_k(\pi)&=&{\rm max}(\mathbb{Z}\cap\{-\langle h_k, \pi(t) \rangle | 0\leq t\leq1 \})\\
&=& -{\rm min}(\mathbb{Z}\cap\{\langle h_k, \pi(t) \rangle | 0\leq t\leq1 \})\\
&=& -{\rm min} \{y_j | j\in[0,r-1]_{\rm int}, \tau_{j+1}\Lambda_i(h_k)>0\}\cup\{y_r\}.
\end{eqnarray*}
For $j\in[1,r]_{\rm int}$, we put
\[
y'_j:=\sum^r_{l=j} \langle h_k,\tau_l\Lambda_i \rangle(a_l-a_{l-1})=y_r-y_{j-1}.
\]

On the other hand, by (\ref{pr-0}), (\ref{pr-00}), the definition of monomial realizations in Sect.\ref{monosec} and
\[
\Phi(\pi)=
\prod^{r}_{l=1}
\left(M(v_{\tau_l\Lambda_i}) \right)^{a_l-a_{l-1}},
\]
we see that
\[
\varphi_k(\Phi(\pi))
={\rm max}(\{ y'_l | l\in[1,r]_{\rm int},\ \tau_l\Lambda_i(h_k)>0 \}\cup\{0\}).
\]
Therefore,
\begin{eqnarray*}
& &\varepsilon_k(\Phi(\pi))\\
&=&
\varphi_k(\Phi(\pi))-
\langle h_k,{\rm wt}(\Phi(\pi))\rangle\\
&=&{\rm max}(\{ y'_l | l\in[1,r]_{\rm int},\ \tau_l\Lambda_i(h_k)>0 \}\cup\{0\})-y_r\\
&=&{\rm max}(\{ -y_{l-1} | l\in[1,r]_{\rm int},\ \tau_l\Lambda_i(h_k)>0 \}\cup\{-y_r\}) \\
&=&{\rm max}(\{ -y_{l} | l\in[0,r-1]_{\rm int},\ \tau_{l+1}\Lambda_i(h_k)>0 \}\cup\{-y_r\}) \\
&=& \varepsilon_k(\pi).
\end{eqnarray*}

\underline{${\rm wt}(\Phi(\pi))={\rm wt}(\pi)$ in $P_{\rm cl}$}

\vspace{2mm}

We fix $k\in\{1,2\}$ and take $\pi=(\tau_1,\cdots,\tau_r ; a_0=0,a_1,\cdots,a_{r-1},a_r=1)\in\mathbb{B}(\Lambda_i)$. Then
\[
{\rm wt}(\pi)=
\sum^{r}_{l=1} (a_l-a_{l-1})\tau_l\Lambda_i .
\]
Since
\[
\Phi(\pi)=
\prod^{r}_{l=1}
\left(M(v_{\tau_l\Lambda_i}) \right)^{a_l-a_{l-1}}
\]
and $M(v_{\tau_l\Lambda_i})$ has the weight $\tau_l\Lambda_i\in P_{\rm cl}$, we get
\[
{\rm wt}(\Phi(\pi))=\sum^{r}_{l=1} (a_l-a_{l-1})\tau_l\Lambda_i
={\rm wt}(\pi).
\]

\underline{$\tilde{f}_k(\Phi(\pi))=\Phi(f_k\pi)$}

\vspace{2mm}

Let $\pi=(\tau_1,\cdots,\tau_r ; a_0=0,a_1,\cdots,a_{r-1},a_r=1)\in\mathbb{B}(\Lambda_i)$.
We consider the action of $f_k$ and use the same notation as in subsection \ref{3-2}:
\[
Q:={\rm min}(\mathbb{Z}\cap\{\langle h_k, \pi(t) \rangle | 0\leq t\leq1 \})=-\varepsilon_k(\pi),
\]
\[
p:={\rm max}\{l\in[0,r]_{\rm int} | y_l=Q\}.
\]
We assume that
\[
\langle h_k,{\rm wt}(\pi)\rangle>Q=-\varepsilon_k(\pi)
\]
so that $\varphi_k(\Phi(\pi))>0$. In particular, it holds $p<r$.
Now, we know
$y_{p+2},y_{p+4},y_{p+6},\cdots\in\mathbb{Z}_{\geq Q+1}$ and $y_{p+1}>y_{p+2}$, $y_{p+3}>y_{p+4}, y_{p+2}$, $y_{p+5}>y_{p+6},y_{p+4}$ (Theorem \ref{Sagthm}). Note that if $p>0$ then $\langle h_k,\tau_p\Lambda_i\rangle<0$ and $s_k\tau_p=\tau_{p+1}$.
Therefore, we see that
\[
x:={\rm min}\{l\in[p,r]_{\rm int} | \langle h_k, \pi(t) \rangle\geq Q+1\ {\rm if}\ t\geq a_l\}=p+1.
\]

We consider a division into cases:

\vspace{3mm}

\underline{Case 1. $p=0$, $h_k(a_1)=y_1=Q+1$}

\vspace{5mm}

\underline{Case 2. $p=0$, $h_k(a_1)=y_1>Q+1$}

\vspace{5mm}

\underline{Case 3. $p>0$, $h_k(a_{p+1})=y_{p+1}=Q+1$}

\vspace{5mm}

\underline{Case 4. $p>0$, $h_k(a_{p+1})=y_{p+1}>Q+1$}

\vspace{4mm}

We define integers $P_{p+1}$ and $P_{p}$ as
\[
\tau_{p+1}\Lambda_i=P_{p+1}\Lambda_k-P_{p}\Lambda_{k'}.
\]
Since $\langle h_k,\tau_{p+1}\Lambda_i\rangle>0$, we obtain $P_{p+1}>0$, $P_p\geq0$. 

\underline{Case 1. $p=0$, $h_k(a_1)=y_1=Q+1$}

\vspace{3mm}

If $r>1$ then $Q+1=h_k(a_1)>h_k(a_2)\geq Q+1$ by $y_{p+1}>y_{p+2}$, which is a contradiction.
Thus, it holds $r=1$, that is, $\pi=(\tau_1;0,1)$. We see $Q=y_0=0$ and $y_1=\langle h_k,\tau_1\Lambda_i \rangle=1$ so that
Definition \ref{opdef} yields
\begin{equation}\label{case1-1}
f_k(\pi)=(s_k\tau_1;0,1).
\end{equation}

\vspace{2mm}

\underline{Case 2. $p=0$, $h_k(a_1)=y_1>Q+1$}

\vspace{3mm}

If $a_0=0<a<a_1$, we have
\[
h_k(a)=h_k(a_0)+ \langle h_k,\tau_1\Lambda_i\rangle(a-a_0)
= Q+P_{1}a.
\]
Thus, the equation $h_k(a)=Q+1$ in Definition \ref{opdef}
means $a=\frac{1}{P_1}$, which yields
\begin{equation}\label{case1-2}
f_k(\pi)=(s_k\tau_1,\tau_1,\tau_2,\cdots,\tau_r;0,\frac{1}{P_1},a_1,a_2,\cdots,a_{r-1},1).
\end{equation}

\underline{Case 3. $p>0$, $h_k(a_{p+1})=y_{p+1}=Q+1$}

\vspace{3mm}

In this case, if $p\leq r-2$ then $Q+1=h_k(a_{p+1})>h_k(a_{p+2})\geq Q+1$, which is a contradiction.
Thus, $p=r-1$ and
\begin{equation}\label{case1-3}
f_k(\pi)=(\tau_1,\cdots,\tau_{r-1};0,a_1,\cdots,a_{r-2},1).
\end{equation}

\vspace{2mm}

\underline{Case 4. $p>0$, $h_k(a_{p+1})=y_{p+1}>Q+1$}

\vspace{3mm}

Just as in Case 2, if $a_p<a<a_{p+1}$ then by the definition of $h_k$ it holds
\[
h_k(a)=h_k(a_p)+ \langle h_k,\tau_{p+1}\Lambda_i\rangle(a-a_p)
= Q+P_{p+1}(a-a_p).
\]
Hence, the equation $h_k(a)=Q+1$ means $a=a_p+\frac{1}{P_{p+1}}$ and
\begin{equation}\label{case1-4}
f_k(\pi)=(\tau_1,\tau_2,\cdots,\tau_r;0,a_1,\cdots,a_{p-1}, a_p+\frac{1}{P_{p+1}},a_{p+1},\cdots,a_{r-1},1).
\end{equation}

Next, let us consider the action of $\tilde{f}_k$ on $\Phi(\pi)$.
Recall that
\begin{eqnarray}
p&=&{\rm max}\{l\in[0,r]_{\rm int} | y_l=Q \}\nonumber\\
&=&{\rm max}\{l\in[0,r]_{\rm int} | y_r-y_l=y_r-Q(=\varphi_k(\pi)) \}.\label{pr-3}
\end{eqnarray}
Note that $y_r-y_l$ is equal to $ \langle h_k, {\rm wt}(\prod^{r}_{j=l+1}M(v_{\tau_j\Lambda_i})^{a_j-a_{j-1}}) \rangle$.
Let $n_{f}$ be the left index of numerator of $M(v_{\tau_{p+1}\Lambda_i})$, that is,
\[
M(v_{\tau_{p+1}\Lambda_i})=\frac{X^{P_{p+1}}_{n_f,k}}{X^{P_{p}}_{n',k'}}
\]
with some $n'\in\mathbb{Z}$.
Since we have shown $\varphi_k(\pi)=\varphi_k(\Phi(\pi))$,
by the definition of monomial realizations and (\ref{pr-3}), it follows
\[
\tilde{f}_k\Phi(\pi)=\Phi(\pi)\cdot A_{n_f,k}^{-1}=
\Phi(\pi)\cdot \frac{X^{-a_{k',k}}_{n_f+\delta_{2,k},k'}}{X_{n_f,k}X_{n_f+1,k}}.
\]

\underline{Case 1. $p=0$, $h_k(a_1)=y_1=Q+1$}

\vspace{2mm}

In this case, $\pi=(\tau_1;0,1)$, $Q=0$ and $\langle h_k,\tau_1\Lambda_i\rangle=h_k(a_1)=1$. Thus, $P_1=1$ so that
$M(v_{\tau_1\Lambda_i})=\frac{X_{n_f,k}}{X^{P_0}_{n',k'}}$.
It follows from (\ref{case1-1}) that
\[
\Phi(f_k\pi)=M(v_{s_k\tau_1\Lambda_i}),\quad
\tilde{f}_k\Phi(\pi)=\tilde{f}_kM(v_{\tau_1\Lambda_i})=M(v_{s_k\tau_1\Lambda_i}).
\]
Therefore, we get $\Phi(f_k\pi)=\tilde{f}_k\Phi(\pi)$.

\vspace{2mm}

\underline{Case 2. $p=0$, $h_k(a_1)=y_1>Q+1$}

\vspace{3mm}

We get
\begin{equation}\label{pr-4}
M(v_{\tau_1\Lambda_i})=\frac{X^{P_1}_{n_f,k}}{X^{P_0}_{n',k'}}
\end{equation}
and
\begin{equation}\label{pr-5}
M(v_{s_k\tau_1\Lambda_i})=M(v_{\tau_1\Lambda_i})A^{-P_1}_{n_f,k}.
\end{equation}
Combining with (\ref{case1-2}),
\[
\Phi(\pi)=M(v_{\tau_1\Lambda_i})^{a_1-a_0}
M(v_{\tau_2\Lambda_i})^{a_2-a_1}\cdots M(v_{\tau_r\Lambda_i})^{a_r-a_{r-1}},
\]
(\ref{pr-4}) and (\ref{pr-5}), we have
\begin{eqnarray*}
\Phi(f_k\pi)
&=& M(v_{s_k\tau_1\Lambda_i})^{\frac{1}{P_1}}M(v_{\tau_1\Lambda_i})^{a_1-\frac{1}{P_1}}
M(v_{\tau_2\Lambda_i})^{a_2-a_1}\cdots M(v_{\tau_r\Lambda_i})^{a_r-a_{r-1}}\\
&=& M(v_{s_k\tau_1\Lambda_i})^{\frac{1}{P_1}}M(v_{\tau_1\Lambda_i})^{-\frac{1}{P_1}} \Phi(\pi)\\
&=& (M(v_{\tau_1\Lambda_i})A^{-P_1}_{n_f,k})^{\frac{1}{P_1}}M(v_{\tau_1\Lambda_i})^{-\frac{1}{P_1}} \Phi(\pi)\\
&=& \Phi(\pi)A^{-1}_{n_f,k}=\tilde{f}_k\Phi(\pi).
\end{eqnarray*}

\underline{Case 3. $p>0$, $h_k(a_{p+1})=y_{p+1}=Q+1$}

\vspace{3mm}

In this case, the previous discussion says $p=r-1$. 
The definition of $p$ implies $y_p=y_{r-1}=Q$.
By the assumption, $y_r=Q+1=y_{r-1}+1$ so that 
\begin{equation}\label{pr-6}
1=y_r-y_{r-1}=(1-a_{r-1})\langle h_k,\tau_r\Lambda_i\rangle=(1-a_{r-1})P_r.
\end{equation}
Since $M(v_{\tau_r\Lambda_i})=\frac{X^{P_{r}}_{n_f,k}}{X^{P_{r-1}}_{n',k'}}$, it holds
\begin{equation}\label{pr-7}
M(v_{\tau_r\Lambda_i})A^{-P_{r}}_{n_f,k}
=M(v_{s_k\tau_r\Lambda_i})=M(v_{\tau_{r-1}\Lambda_i}).
\end{equation}
It follows by (\ref{case1-3}),
\[
\Phi(\pi)=
M(v_{\tau_1\Lambda_i})^{a_1-a_0}
\cdots M(v_{\tau_{r-2}\Lambda_i})^{a_{r-2}-a_{r-3}} M(v_{\tau_{r-1}\Lambda_i})^{a_{r-1}-a_{r-2}}
M(v_{\tau_r\Lambda_i})^{1-a_{r-1}},
\]
(\ref{pr-6}) and (\ref{pr-7}) that
\begin{eqnarray*}
\Phi(f_k\pi)
&=& 
M(v_{\tau_1\Lambda_i})^{a_1-a_0}
\cdots M(v_{\tau_{r-2}\Lambda_i})^{a_{r-2}-a_{r-3}}
M(v_{\tau_{r-1}\Lambda_i})^{1-a_{r-2}}\\
&=& \Phi(\pi) M(v_{\tau_{r}\Lambda_i})^{a_{r-1}-1}
M(v_{\tau_{r-1}\Lambda_i})^{1-a_{r-1}}\\
&=&\Phi(\pi) M(v_{\tau_{r}\Lambda_i})^{a_{r-1}-1}
(M(v_{\tau_r\Lambda_i})A^{-P_{r}}_{n_f,k})^{1-a_{r-1}}\\
&=&\Phi(\pi) A^{-P_{r}(1-a_{r-1})}_{n_f,k}=\Phi(\pi) A^{-1}_{n_f,k}
=\tilde{f}_k\Phi(\pi).
\end{eqnarray*}

\underline{Case 4. $p>0$, $h_k(a_{p+1})=y_{p+1}>Q+1$}

\vspace{3mm}

In this case, we get
\[
M(v_{\tau_{p+1}\Lambda_i})=\frac{X^{P_{p+1}}_{n_f,k}}{X^{P_p}_{n',k'}},
\]
\begin{equation}\label{pr-8}
M(v_{\tau_p\Lambda_i})=
M(v_{s_k\tau_{p+1}\Lambda_i})=\tilde{f}^{P_{p+1}}_kM(v_{\tau_{p+1}\Lambda_i})=M(v_{\tau_{p+1}\Lambda_i})A^{-P_{p+1}}_{n_f,k}.
\end{equation}
Using (\ref{case1-4}) and (\ref{pr-8}), we have
\begin{eqnarray*}
\Phi(f_k(\pi))&=& 
M(v_{\tau_1\Lambda_i})^{a_1-a_0}\cdots M(v_{\tau_{p-1}\Lambda_i})^{a_{p-1}-a_{p-2}}
\\
& &M(v_{\tau_{p}\Lambda_i})^{a_{p}+\frac{1}{P_{p+1}}-a_{p-1}}
M(v_{\tau_{p+1}\Lambda_i})^{a_{p+1}-a_{p}-\frac{1}{P_{p+1}}}\\
& &M(v_{\tau_{p+2}\Lambda_i})^{a_{p+2}-a_{p+1}}
\cdots M(v_{\tau_{r}\Lambda_i})^{a_{r}-a_{r-1}}\\
&=& \Phi(\pi) M(v_{\tau_p\Lambda_i})^{\frac{1}{P_{p+1}}}M(v_{\tau_{p+1}\Lambda_i})^{-\frac{1}{P_{p+1}}}\\
&=& \Phi(\pi) (M(v_{\tau_{p+1}\Lambda_i})A^{-P_{p+1}}_{n_f,k})^{\frac{1}{P_{p+1}}} M(v_{\tau_{p+1}\Lambda_i})^{-\frac{1}{P_{p+1}}}\\
&=& \Phi(\pi)A^{-1}_{n_f,k}=\tilde{f}_k\Phi(\pi).
\end{eqnarray*}

By $\Phi(e;0,1)=M(v_{\Lambda_i})=X_{s,i}\in \mathcal{M}(\Lambda_i)$ and above argument, we see that
\[
\Phi(\mathbb{B}(\Lambda_i))\subset \mathcal{M}(\Lambda_i).
\]
Let us prove $\Phi$ is injective.
Taking
\[
\pi=(\tau_1,\cdots,\tau_r;0,a_1,\cdots,a_{r-1},1)\in
\mathbb{B}(\Lambda_i),
\]
\[
\pi'=(\tau'_1,\cdots,\tau'_m;0,a'_1,\cdots,a'_{m-1},1)\in
\mathbb{B}(\Lambda_i)
\]
and we suppose that $\Phi(\pi)=\Phi(\pi')$ and prove $\pi=\pi'$. One may assume $\tau_1\neq e$, $\tau'_1\neq e$.
Let us focus on the denominators of $M(v_{\tau_1\Lambda_i})$ and $M(v_{\tau'_1\Lambda_i})$ and put them as
\[
M(v_{\tau_1\Lambda_i})=\frac{X^{a}_{b,c}}{X^{P_{t,l}}_{t,l}},\quad
M(v_{\tau'_1\Lambda_i})=\frac{X^{a'}_{b',c'}}{X^{P_{t',l'}}_{t',l'}},
\]
with some $t,t'\in\mathbb{Z}$, $l,l'\in\{1,2\}$, $P_{t,l}, P_{t',l'}\in\mathbb{Z}_{>0}$
and some monomials $X^{a}_{b,c}$ and $X^{a'}_{b',c'}$.
Considering the explicit forms of $M(v_{\tau\Lambda_i})$ in (\ref{pr-0}), $(\ref{pr-00})$,
we see that $(t,l)$ (resp. $(t',l')$) is the maximum double index appearing in $\Phi(\pi)$ (resp. $\Phi(\pi')$)
in the order 
\[
(s,1)<(s,2)<(s+1,1)<(s+1,2)<(s+2,1)<(s+2,2)<\cdots.
\]
Since we supposed $\Phi(\pi)=\Phi(\pi')$, it holds $(t,l)=(t',l')$,
$M(v_{\tau_1\Lambda_i})^{a_1}=M(v_{\tau'_1\Lambda_i})^{a_1'}$ so that $\tau_1=\tau'_1$ and $a_1=a_1'$.
By 
\[
(M(v_{\tau_1\Lambda_i})^{a_1})^{-1}\Phi(\pi)=(M(v_{\tau_1\Lambda_i})^{a_1})^{-1}\Phi(\pi')
\]
and a similar argument,
it holds $M(v_{\tau_2\Lambda_i})^{a_2-a_1}=M(v_{\tau'_2\Lambda_i})^{a_2'-a_1'}$ so that $\tau_2=\tau'_2$ and $a_2=a_2'$.
Repeating this argument, we see that $r=m$, $a_j=a'_j$, $\tau_j=\tau'_j$ ($j\in[1,r]_{\rm int}$) and $\pi=\pi'$.
We proved $\Phi$ is injective from $\mathbb{B}(\Lambda_i)$ to $\mathcal{M}(\Lambda_i)$.

For any $\tilde{f}_{i_t}\cdots\tilde{f}_{i_2}\tilde{f}_{i_1}X_{s,i}\in\mathcal{M}(\Lambda_i)$, since the map $\Phi$ preserves
$\varphi_k$, it holds $f_{i_t}\cdots f_{i_2}f_{i_1}(1;0,1)\in \mathbb{B}(\Lambda_i)$ and
\[
\Phi(f_{i_t}\cdots f_{i_2}f_{i_1}(e;0,1))=\tilde{f}_{i_t}\cdots\tilde{f}_{i_2}\tilde{f}_{i_1}X_{s,i}.
\]
Thus $\Phi$ is surjective.

\vspace{3mm}

For $\pi\in \mathbb{B}(\Lambda_i)$, it holds $e_k\pi=0$ if and only if $\tilde{e}_k\Phi(\pi)=0$. 
If $e_k\pi\neq0$ then we can write $e_k\pi=f_{i_t}\cdots f_{i_1}(e;0,1)$ with some $i_1,\cdots,i_t\in I$.
Above argument implies
\[
\Phi(e_k\pi)=\tilde{f}_{i_t}\cdots \tilde{f}_{i_1}X_{s,i},\quad
\tilde{e}_k\Phi(\pi)=\tilde{e}_k\Phi(f_kf_{i_t}\cdots f_{i_1}(e;0,1))
=\tilde{e}_k\tilde{f}_k\tilde{f}_{i_t}\cdots \tilde{f}_{i_1}X_{s,i}
=\tilde{f}_{i_t}\cdots \tilde{f}_{i_1}X_{s,i}
\]
so that the relation $\Phi(e_k \pi)=\tilde{e}_k\Phi(\pi)$ also holds.
Thus $\Phi$ is an isomorphism
between $\mathbb{B}(\Lambda_i)$ and $\mathcal{M}(\Lambda_i)$.\qed


\begin{thebibliography}{9}


\bibitem[JKKS]{JKKS} 
K.Jeong, S.-J.Kang, J.-A.Kim, D.-U.Shin,
Crystals and Nakajima monomials for quantum
generalized Kac–Moody algebras,
J. Algebra, 319, no.9, pp.3732--3751(2008).





\bibitem[Jos]{Jos} A.Joseph, {\it Quantum Groups and Their Primitive Ideals}, Ergebnisse der Mathematik und ihrer Grenzgebiete (3),
29, Springer–Verlag, Berlin (1995).


\bibitem[Kas90]{Kas90} M.Kashiwara, Crystalling the $q$-analogue of universal 
              enveloping algebras, Comm. Math. Phys.,
            {\it 133}, 249--260 (1990).


\bibitem[Kas93]{Kas93} M.Kashiwara, 
The crystal base and Littelmann's refined Demazure character formula,
Duke Math. J., 71, no 3, pp.839--858 (1993).

\bibitem[Kas96]{Kas96} M.Kashiwara, Similarity of crystal bases,
Lie algebras and their representations (Seoul, 1995), 
Contemp. Math., 194,
American Mathematical Society, Providence, RI, pp.177–186 (1996).


\bibitem[Kas03]{Kas03} M.Kashiwara, Realizations of crystals, 
Combinatorial and geometric representation theory (Seoul, 2001), 
Contemp. Math., 325, Amer. Math. Soc., Providence, RI, pp.133--139 (2003).


\bibitem[Lit94]{Lit94} P.Littelmann,
A Littlewood-Richardson rule for symmetrizable Kac-Moody algebras, 
Invent. Math.116, no.1-3, pp.329--346 (1994).


\bibitem[Lit95]{Lit95} P.Littelmann,
Paths and root operators in representation theory, 
Ann. of Math. (2)142, no.3, pp.499--525 (1995).


\bibitem[Nak]{Nak03} H.Nakajima, $t$-analogs of $q$-characters of quantum affine algebras of type $A_n$, $D_n$,
Combinatorial and geometric representation theory (Seoul, 2001),
Contemp. Math., 325, Amer. Math. Soc., Providence, RI, pp.141--160 (2003).

\bibitem[Sag]{Sag} D.Sagaki, [Path model and various related results] path model to soreni kanrenshita shokekka (in Japanese),
RIMS K\={o}ky\^{u}roku, 1183, pp.35-51 (2001).




\bibitem[San]{San} Y.B.Sanderson,
Dimensions of Demazure modules for rank two affine Lie algebras, Compositio Math.101, no.2, pp.115--131 (1996).






\end{thebibliography}
\end{document}